\documentstyle[11pt]{article}

\begin{document}
\title{On the Classification of K3 Surfaces with Nine Cusps}
\author{W. Barth} 
\maketitle 

\begin{center}
Mathematisches Institut der Universit\"{a}t

Bismarckstr. 1 1/2, D 91054 Erlangen
\end{center}

\def\Rn{\R^n}

\def \qed {\hspace*{\fill}\frame{\rule[0pt]{0pt}{8pt}\rule[0pt]{8pt}{0pt}}\par}

\def\R{{\rm I\!R}} 
\def\N{{\rm I\!N}} 
\def\F{{\rm I\!F}}
\def\P{{\rm I\!P}}
\def\E{{\mathchoice {\rm 1\mskip-4mu l} {\rm 1\mskip-4mu l}
{\rm 1\mskip-4.5mu l} {\rm 1\mskip-5mu l}}}
\def\C{{\mathchoice {\setbox0=\hbox{$\displaystyle\rm C$}\hbox{\hbox
to0pt{\kern0.4\wd0\vrule height0.9\ht0\hss}\box0}}
{\setbox0=\hbox{$\textstyle\rm C$}\hbox{\hbox
to0pt{\kern0.4\wd0\vrule height0.9\ht0\hss}\box0}}
{\setbox0=\hbox{$\scriptstyle\rm C$}\hbox{\hbox
to0pt{\kern0.4\wd0\vrule height0.9\ht0\hss}\box0}}
{\setbox0=\hbox{$\scriptscriptstyle\rm C$}\hbox{\hbox
to0pt{\kern0.4\wd0\vrule height0.9\ht0\hss}\box0}}}}
\def\Q{{\mathchoice {\setbox0=\hbox{$\displaystyle\rm Q$}\hbox{\raise
0.15\ht0\hbox to0pt{\kern0.4\wd0\vrule height0.8\ht0\hss}\box0}}
{\setbox0=\hbox{$\textstyle\rm Q$}\hbox{\raise
0.15\ht0\hbox to0pt{\kern0.4\wd0\vrule height0.8\ht0\hss}\box0}}
{\setbox0=\hbox{$\scriptstyle\rm Q$}\hbox{\raise
0.15\ht0\hbox to0pt{\kern0.4\wd0\vrule height0.7\ht0\hss}\box0}}
{\setbox0=\hbox{$\scriptscriptstyle\rm Q$}\hbox{\raise
0.15\ht0\hbox to0pt{\kern0.4\wd0\vrule height0.7\ht0\hss}\box0}}}}
\def\Z{{\mathchoice {\hbox{$\sf\textstyle Z\kern-0.4em Z$}}
{\hbox{$\sf\textstyle Z\kern-0.4em Z$}}
{\hbox{$\sf\scriptstyle Z\kern-0.3em Z$}}
{\hbox{$\sf\scriptscriptstyle Z\kern-0.2em Z$}}}}

\newcommand{\D}{\displaystyle} 
\newcommand{\vv}{\vspace{0.3cm}} 
\newcommand{\vV}{\vspace{0.5cm}} 
\newcommand{\VV}{\vspace{2cm}}

\begin{abstract}
By a $K3$-surface with nine cusps I mean a 
compact complex surface with nine
isolated double points $A_2$, but otherwise smooth, such that
its minimal desingularisation is a $K3$-surface.
In an earlier paper I showd that each such surface is a 
quotient of a complex torus by a cyclic group of order three.
Here I try to classify these $K3$-surfaces, using the period
map for complex tori. In particular I show:

A $K3$-surface with nine cusps carries polarizations only of
degrees $0$ or $2$ modulo $6$. This implies in particular that
there is no quartic surface in projective three-space with
nine cusps. (T. Urabe pointed out to me how to deduce
this from Nikulin's Thm. 1.12.2 in [N2].)

In an appendix I give explicit 
equations of quartic surfaces in three-space
with eight cusps. 
 

\end{abstract}

\tableofcontents

\setcounter{section}{-1}

\VV
\section{Introduction}
In [B] it was shown that each compact complex surface with nine
cusps ($A_2$-type double points), but no further singularities,
such that its minimal desingularization is a $K3$-surface, arises
as a $3:1$ quotient of a complex torus. It was claimed there, that 
there are nonalgebraic surfaces of this kind, and it was suggested that
the examples of [BL] are the only algebraic surfaces of this type: 
double covers of the plane, branched over
the dual sextic to a smooth cubic curve.

The aim of this note is to prove the first claim, and to 
show that there are lots of other examples than those in [BL]. 
In fact, the author was informed by P. Vanhaecke of
his joint note [BV] with J. Bertin, where sextic surfaces in $\P_4$
are constructed, complete intersections of a quartic 
and a cubic hypersurface, which have nine cusps. They
are $3:1$ quotients of Jacobians of genus-2 curves with
an automorphism of order three.

It is not an accident, that the surfaces of [BL] carry a polarization
of degree two, while the surfaces of [BV] have a polarization
of degree six. It is shown below, that only polarizations of
degrees $0$ or $2$ modulo $6$ appear on algebraic $K3$-surfaces
with nine cusps. This has the following consequence
(T. Urabe pointed out to me that it is also a consequence
of Nikulin's theorem 1.12.2 in [N2], although it seems not
to have been noted before explicitely):

{\em There is no quartic surface in $\P_3$ with nine cusps, and no further
singularities (or with other isolated rational double points only).}

This leads to the following obvious question:
Are there quartic surfaces in $\P_3$ with eight cusps?
The existenc follows indeed from Urabe's partial classification
of quartic surfaces with simple singularities [U, thm. 0.2]
But it seems quite hard to write down their equations directly.
Fortunately [BV] describe their sextic surfaces in $\P_4$
with nine cusps quite explicitely. So it is fairly easy to project
them from one of their cusps into $\P_3$. There probably is some general reason
for the fact, that none of the projections degenerates. They all have
eight cusps and no further singularities. However I check this
by simple, but tedious direct computation (see the appendix).   

I use the period ($=$ holomorphic $2$-form) of the covering
two-dimensional torus to classify, at least to some extent, complex
two-dimensional tori admitting a cyclic symmetry group of order three
such that the quotient is a $K3$-surface with nine cusps. The period domain
$\Omega$
is an open dense subset in a smooth quadric of dimension two, in particular
it is connected. The moduli space for pairs $(A,t)$ with $A$ a complex torus
of dimension two and $t$ an automorphism of order three as above,
is a quotient of $\Omega$ by an infinite arithmetic group. The general surface
$A$ of this type is not algebraic.

Fixing a polarization $\alpha$ on $A$ reduces the period domain to
a curve $\Omega_{\alpha} \subset \Omega$. Funny enough this 
curve consists of two disjoint copies of the upper 
half-plane, interchanged by conjugation of the complex
structure on $A$. Here I mean by a polarization
a fixed divisor class $\alpha \in H^2(A,\Z)$ with $\alpha^2>0$. 
Of course there are
infinitely many different classes $\alpha$, which are equivalent under the
group $SL(4,\Z)$. To induce a polarization on the quotient
$X=A/t$, the class $\alpha$ has at least to be $t$-invariant. And the
group identifying isomorphic triplets $A,t,\alpha$ is the subgroup
$G \subset SL(4,\Z)$ of elements commuting with $t$. Unfortunately it seems
quite difficult to classify $G$-orbits on the set of classes $\alpha$
of fixed square $\alpha^2$. At least it is a problem which I
cannot solve.  

Finally I give an example of algebraic tori $A$ with automorphism $t$
which are simple, i.e., not isogenous to a product of elliptic curves.

Just as the note [B] was essentially parallel to the first pages of
Nikulin's paper [N1], 
which treats the natural involution $a \mapsto -a$,
the basic method here is parallel to Remark 2 in section 1 of [N1].
Using the discriminant of its quadratic form we identify
the orthogonal complement in the $K3$-lattice of the sub-lattice
$I$ spanned by the 18 classes of the rational curves resolving
the nine cusps.  

\VV
\section{Topology}

\subsection{The action on $H_1(A,\Z)$}

In this section let $A=\C^2/\Gamma$ be a complex torus
of dimension two. Let
$t:A \to A$ be an automorphism of order three having the
origin as an isolated fixed point. 
Asume also that $t^*(\alpha)=\alpha$, where
$\alpha$ is the holomorphic $2$-form induced by the constant
form $dz_1 \wedge dz_2$ on $\C^2$.

$t$ induces a linear automorphism $\tilde{t}:\C^2 \to \C^2$
on the universal covering, with determinant $=1$. Since the origin
is an isolated fixed point, it has no eigenvalue $=1$. 
So it must have the two eigenvalues $\omega=e^{2\pi i/3}, \omega^2$.

In this section I want to identify the action of $t$ on
the homology $H_1(A,\Z)$, or what is the same, the action
of $\tilde{t}$ on the lattice $\Gamma \simeq \Z^4 \subset \C^2$.

As there is no $t$-invariant real line in $\C^2$, 
nor a $t$-invariant real subvector space of dimension three, there
is no invariant sub-lattice in $H_1(A,\Z)$ of rank one or three.
But the $t$-orbit of each period $\gamma \in \Gamma$ 
spans an invariant sub-lattice of rank $\leq 3$. This implies
that each orbit $0 \neq \gamma, t(\gamma), t^2(\gamma)$ spans a
sub-lattice of rank two. 

Consider some primitve vector $\alpha \in \Gamma$ and denote
the primitive vector $t(\alpha)$ by $\beta$. Then
$$t(\beta)= p \cdot \alpha + q \cdot \beta \quad \mbox{with} \quad
p,q \in \Q$$
and
$$\left( \begin{array}{cc} 0 & p \\ 1 & q \\ \end{array} \right)^3 =
\left( \begin{array}{cc}
       pq & p(p+q^2) \\ p+q^2 & q(2p+q^2) \\ \end{array} \right)
=\left( \begin{array}{cc} 1 & 0 \\ 0 & 1 \\ \end{array} \right).$$
Here $pq=1$ and $-p=q^2=1/p^2$ imply $p=q=-1$. So the action
of $t$ on the invariant sub-lattice generated by $\alpha$ and $\beta$
is given by the matrix 
$$ \left( \begin{array}{cc} 0 & -1 \\ 1 & -1 \\ \end{array} \right).$$

{\bf Proposition:} {\em There is a $\Z$-basis of the lattice $H_1(A,\Z)$
$$\alpha_1, \; \beta_1,\; \alpha_2, \; \beta_2,$$
in which the action of $t$ is}
$$t(\alpha_1)=\beta_1, \; t(\beta_1)=-\alpha_1-\beta_1,  \quad 
  t(\alpha_2)=\beta_2, \; t(\beta_2)=-\alpha_2-\beta_2.$$

Proof. Consider $H_1(A,\Z)$ as a sub-lattice of the real vector
space $H_1(A,\R)\simeq \R^4$ and choose on this vector space
some $t$-invariant inner product $(-,-)$. Let 
$\alpha_1 \in H_1(A,\Z)$ be some lattice vector of smallest length 
$\parallel \alpha_1 \parallel =\sqrt{(\alpha_1,\alpha_1)} \not= 0$ 
and put $\beta_1:=t(\alpha_1)$. On the plane
spanned by $\alpha_1$ and $\beta_1$ the automorphism $t$ is an isometry
of order three. This implies
$$(\alpha_1,\beta_1) = -\frac{1}{2} \parallel \alpha_1 \parallel^2.$$
If the vectors $\alpha_1$ and $\beta_1$ would not span a primitive
sub-lattice of $H_1(A,\Z)$, there would be some nonzero
lattice vector $\gamma =u \alpha_1+v\beta_1 \in H_1(A,\Z)$ with
$0 \leq u,v <1$. But such a vector would have squared length 
$$(u^2+v^2)\parallel \alpha_1 \parallel^2 +2uv \cdot (\alpha_1, \beta_1) 
=(u^2+v^2-uv) \parallel \alpha_1 \parallel^2.$$
Either $u^2 \leq uv$ or $v^2 \leq uv$, hence
$$u^2+v^2-uv \leq max\{u^2,v^2\} < 1.$$
Such a vector $\gamma \not=0 $ of length $<\parallel \alpha_1 \parallel$
cannot exist. So $\alpha_1$ and $\beta_1$ indeed span a primitive
sub-lattice in $H_1(A,\Z)$.

Now choose a lattice vector $\alpha_2 \in H_1(A,\Z)$ of smallest
distance $\not=0$ from the subvector space generated 
by $\alpha_1$ and $\beta_1$
in $H_1(A,\R)$ with $\beta_2=t(\alpha_2)$. Exactly the same argument
shows that the residues $\bar{\alpha_2}$ and $\bar{\beta_2}$
form a $\Z$-basis of the quotient lattice 
$H_1(A,\Z)/(\Z \cdot \alpha_1 + \Z \cdot \beta_1)$. So 
$\alpha_1,\beta_1,\alpha_2,\beta_2$ form a $\Z$-basis of $H_1(A,\Z)$. \qed

\vv
With respect to a $\Z$-basis as in the proposition the action of
$t$ is given by the matrix
$$T:= \left( \begin{array}{rrrr}
             0 & -1 & 0 &  0 \\
             1 & -1 & 0 &  0 \\
             0 &  0 & 0 & -1 \\
             0 &  0 & 1 & -1 \\  \end{array} \right).$$

\vv
As a sub-lattice of $\C^2$ the lattice $\Gamma$ carries a natural
orientation.

{\bf Proposition:} {\em Each basis 
$\alpha_1,t(\alpha_1),\alpha_2,t(\alpha_2)$ as
above is negatively oriented with respect to the natural orientation
on $\Gamma$.}

Proof. As a complex linear map, $t$ has the two
eigenvalues $\omega$ and $\omega^2$. Let $c_1, c_2 \in \C^2$ be a 
complex basis of eigenvectors, so $t(c_1)=\omega c_1$ and
$t(c_2)=\omega^2 c_2$. Clearly $c_1$ and $\omega c_1$ represent
the natural orientation of the complex line containing these
two vectors, while $c_2$ and $\omega^2 c_2$ represent the opposite
of the natural orientation of their line. This shows:

There is some real basis $\alpha_1=c_1, t(\alpha_1), \alpha_2=c_2,
t(\alpha_2)$ of $\C^2$ representing the opposite of the natural
orientation. The assertion follows, if we prove the next

{\bf Lemma:} {\em All $\R$-bases of the form 
$\alpha_1,t(\alpha_1),\alpha_2,t(\alpha_2)$
of $\C^2$ represent the same orientation.}

Proof. Let us change from the basis 
$\alpha_1,\beta_1=t(\alpha_1),\alpha_2,\beta_2=t(\alpha_2)$ 
to another $\R$-basis 
$\alpha_1',\beta_1'=t(\alpha_1'),\alpha_2',\beta_2'=t(\alpha_2')$ of $\C^2$.

a) If $\alpha_1'$ and $\beta_1'$ span the same real plane as
$\alpha_1$ and $\beta_1$, and if $\alpha_2'$ and $\beta_2'$ span the same
real plane as $\alpha_2$ and $\beta_2$, then the pairs 
$\alpha_i'$ and $\beta_i'$
are obtained from $\alpha_i$ and $\beta_i, \, i=1,2$ in their respective planes
by orientation preserving rotations. So 
$\alpha_1',\beta_1',\alpha_2',\beta_2'$
represent the same orientation as
$\alpha_1,\beta_1,\alpha_2,\beta_2$.

b) If $\R \alpha_1' + \R \beta_1' = \R \alpha_1+\R \beta_1$, but
$\R \alpha_2' + \R \beta_2' \not= \R \alpha_2+ \R \beta_2$, then, by a), 
$\alpha_1',\beta_1',\alpha_2',\beta_2'$ has the same orientation as
$\alpha_1,\beta_1,\alpha_2',\beta_2'$. Write
$$\alpha_2' = a_1 \alpha_1+a_2 \beta_1+a_3 \alpha_2 + a_4 \beta_2.$$
Then
$$\beta_2' = 
-a_2 \alpha_1 +(a_1-a_2) \beta_1 - a_4 \alpha_2 +(a_3-a_4) \beta_2$$
and
$$det(\alpha_1,\beta_1,\alpha_2',\beta_2') = det \left( \begin{array}{cccc}
1 & 0 & a_1 & -a_2  \\ 0 & 1 & a_2 & a_1-a_2   \\
0 & 0 & a_3 & -a_4  \\ 0 & 0 & a_4 & a_3-a_4 \\ \end{array} \right)
= a_3^2+a_4^2-a_3a_4 \geq 0.$$
Again the orientation of 
$\alpha_1',\beta_1',\alpha_2',\beta_2'$ is the same as the one of
$\alpha_1,\beta_1,\alpha_2,\beta_2$.

c) If $\R \alpha_i' + \R \beta_i' \not= \R \alpha_j+\R \beta_j$ for
$i,j=1,2$, then by b) the orientation does not change, if we pass
from $\alpha_1,\beta_1,\alpha_2,\beta_2$ to 
$\alpha_1,\beta_1,\alpha_2',\beta_2'$, and 
if we pass from $\alpha_1,\beta_2,\alpha_2',\beta_2'$ to
$\alpha_1',\beta_1',\alpha_2',\beta_2'$. \qed
 
\vV
\subsection{The action on $H_2(A,\Z)=\bigwedge^2 H_1(A,\Z)$}

On the exterior products of the basis vectors from the last
section the automorphism $t$ acts as follows:
\begin{eqnarray*}
\alpha_1 \wedge \beta_1 &\mapsto & \beta_1 \wedge (-\alpha_1-\beta_1) 
                   = \alpha_1 \wedge \beta_1 \\
\alpha_2 \wedge \beta_2 & \mapsto & \alpha_2 \wedge \beta_2 \\
\alpha_1 \wedge \alpha_2 & \mapsto & \beta_1 \wedge \beta_2 \\
\alpha_1 \wedge \beta_2 &\mapsto & -\beta_1 \wedge \alpha_2 
                                     - \beta_1 \wedge \beta_2 \\
\beta_1 \wedge \alpha_2 &\mapsto & -\alpha_1 \wedge \beta_2 
                                     - \beta_1 \wedge \beta_2 \\
\beta_1 \wedge \beta_2 & \mapsto & \alpha_1 \wedge \alpha_2 
                                   + \alpha_1 \wedge \beta_2 +
                     \beta_1 \wedge \alpha_2 + \beta_1 \wedge \beta_2.
\end{eqnarray*}

Using this table one finds the following $t$-invariant classes
in $H_2(A,\Z)$:
\begin{eqnarray*}
\gamma_1 &:=& -\alpha_1 \wedge \beta_1, \\
\gamma_2 &:=& \alpha_2 \wedge \beta_2, \\
\gamma_3 &:=& \alpha_1 \wedge \beta_2 - \beta_1 \wedge \alpha_2, \\
\gamma_4 &:=& \alpha_1 \wedge \alpha_2+\alpha_1 \wedge \beta_2
                                  +\beta_1 \wedge \beta_2.
\end{eqnarray*} 

The wedge product induces on $H_2(A,\Z)$ an integral, unimodular quadratic form
$$(\alpha, \alpha'):=\frac{\alpha \wedge \alpha'}
                     {\alpha_1 \wedge \alpha_2 \wedge \beta_1 \wedge \beta_2}$$
with discriminant $-1$. (Recall from the last section that
$\alpha_1 \wedge \alpha_2 \wedge \beta_1 \wedge \beta_2$ represents the
natural orientation.)
The matrix $(\gamma_i,\gamma_j)$ is
$$\left( \begin{array}{rrrr}
          0 &  1 & 0 & 0 \\
          1 &  0 & 0 & 0 \\
          0 &  0 & 2 & 1 \\
          0 &  0 & 1 & 2 \\ \end{array} \right)$$
with determinant $-3$. This shows that the invariant classes
$\gamma_1,...,\gamma_4$ span a primitive sub-lattice 
$L_A \subset H_2(A,\Z)$ of rank four.

{\bf Proposition:} {\em The lattice $L_A$ is the sublattice
$H_2(A,\Z)^{inv} \subset H_2(A,\Z)$
of all $t$-invariant classes.}

Proof. As the sub-lattice $L_A \subset H_2(A,\Z)$ is primitive,
it suffices to show that the invariant subspace 
$H_2(A,\R)^{inv} \subset H_2(A,\R)$ has dimension at most four. But from
the table above, exhibiting the action of $t$ on 
the exterior products of the $\alpha_i$ and  $\beta_j$ one reads off that the
action of $t$ on $H_2(A,\Z)$ has trace $=3$. So $H_2(A,\R)^{inv}$
is a proper subspace of $H_2(A,\R)$. Its dimension cannot be five,
since then there would be a $t$-invariant one-dimensional
complement, and all of $H_2(A,\R)$ would be $t$-invariant.
So its dimension is at most four. \qed

\vv
The dual lattice $\check{L}_A\subset H_2(A,\R)$ consists of
all classes $\gamma \in \R \cdot L_A \subset H_2(A,\R)$ with
$$(\alpha,\gamma) \in \Z \mbox{ for all } \alpha \in L_A.$$

{\bf Proposition:} {\em The dual lattice $\check{L}_A$ has a
$\Z$-basis consisting of}
$$\gamma_1, \, \gamma_2, \gamma_3, \frac{1}{3}(\gamma_3+\gamma_4).$$

Proof. Let
$$\gamma = c_1 \gamma_1+...+c_4 \gamma_4 \in \R \cdot L_A, \;
c_i \in \R.$$
Then
$$\begin{array}{lllllll}
(\gamma,\gamma_1) &=& c_2, & \qquad &  (\gamma,\gamma_2) &=&  c_1, \\
(\gamma,\gamma_3) &=& 2c_3+c_4, & & (\gamma,\gamma_4) &=& c_3+2c_4. \\
\end{array}$$
So $\gamma \in \check{L}_A$ if and only if $c_1,c_2 \in \Z$ and if
$$2c_3+c_4 \in \Z, \; c_3+2c_4 \in \Z.$$
But the latter is equivalent with
$$ c_3-c_4 \in \Z, \quad 3(c_3+c_4) \in \Z.$$
\qed

In the basis $\gamma_1,\gamma_2,\gamma_3,(\gamma_3+\gamma_4)/3$ 
the quadratic form has the matrix
$$ \left( \begin{array}{rrrr}
0 & 1 & 0 & 0 \\ 1 & 0 & 0 & 0 \\ 0 & 0 & 2 & 1 \\ 0 & 0 & 1 & 2/3 \\
\end{array} \right)$$
and the discriminant $-1/3$. 

\vV
\subsection{The maps $q_*$
and $q^*$}

As in [B], let $q:A \to A/t=X$ be the quotient map. Let $\tilde{A}$
be the blow-up of $A$ in the nine fixed points of $t$ and 
$\tilde{X}$ the minimal desingularization of $X$. One has
$H_2(\tilde{A},\Z) = \Z^9 \perp H_2(A,\Z)$ with $\Z^9$
spanned by the classes of the nine exceptional curves.

As in [B] denote by $I \subset H_2(\tilde{X},\Z)$ the lattice spanned by
the eighteen $(-1)$-curves resolving the nine cusps on $X$.
It is a lattice of rank $18$ and discriminant $3^9$. In [B] it was
shown that it is contained in a primitive sublattice
$\bar{I} \subset H_2(\tilde{X},\Z)$ with
the quotient $\bar{I}/I$ being a group of order $3^3$. This
implies that the discriminant $d(\bar{I})$ of the quadratic form on
$\bar{I}$ is $3^9/(3^3)^2=3^3$.

The orthogonal complement 
$L_X:=I^{\perp}=\bar{I}^{\perp} \subset H_2(\tilde{X},\Z)$
is a lattice of rank $22-18=4$. Since $H_2(\tilde{X},\Z)$
is unimodular with discriminant $-1$ the discriminant
of $L_X$ is
$$d(L_X) = -d(\bar{I}) = -3^3.$$   
The induced map $\tilde{q}_*$ clearly maps $H_2(A,\Z)$  
into $L_X$.

We identify
$$H_2(\tilde{A},\Z)=H^2(\tilde{A},\Z), \quad
H_2(\tilde{X},\Z)=H^2(\tilde{X},\Z)$$
via Poincare-duality. Then we get a morphism
$$\tilde{q}^*:H_2(\tilde{X},\Z) \to H_2(\tilde{A},\Z).$$
It satisfies
$$(\tilde{q}^* \xi, \tilde{q}^* \xi')= 
3 \cdot (\xi, \xi') \mbox{ for all } 
\xi, \xi' \in H_2(\tilde{X},\Z),$$
and the projection formula
$$ (\tilde{q}_* \alpha, \xi) = (\alpha, \tilde{q}^* \xi)
\mbox{ for all } \alpha \in H_2(\tilde{A},\Z), \,
\xi \in H_2(\tilde{X},\Z).$$

{\bf Proposition:} {\em The map $\tilde{q}_* : H_2(\tilde{A}, \Z)
\to H_2(\tilde{X},\Z)$ induces an isomorphism
$$\check{L}_A(3) \to L_X$$
and the map $\tilde{q}^*:H_2(\tilde{X},\Z) \to H_2(\tilde{A},\Z)$
induces an isomorphism}
$$L_X(3) \to 3 \cdot \check{L}_A.$$

Proof. Since all classes $q^*(\xi), \, \xi \in H_2(\tilde{X},\Z)$, 
are $t$-invariant, $\tilde{q}^*(L_X)$ 
is a sub-lattice of $L_A$. The endomorphism
$\tilde{q}_* \, \tilde{q}^* : L_X \to L_X$ has the property
$$\tilde{q}_*\tilde{q}^*(\xi)=3\xi$$
for all $\xi \in L_X$.
This implies that the maps $\tilde{q}^*:L_X \to L_A$ and 
$\tilde{q}_*:L_A \to L_X$
are injective. So $q_*(L_A) \subset L_X$
is a sublattice of rank four, the rank of $L_A$, 
and $\tilde{q}^*(L_X)$ spans $L_A$ over $\Q$.

For all $\xi, \xi' \in L_X$ we have
$$(\tilde{q}^*\xi,\tilde{q}^*\xi') = 3 (\xi,\xi')
\mbox{ and } 
(\tilde{q}_*\tilde{q}^*\xi,\tilde{q}_*\tilde{q}^*\xi') = 
(3\xi,3\xi')=3(\tilde{q}^*\xi,\tilde{q}^*\xi').$$
This implies
$$(\tilde{q}_*\alpha, \tilde{q}_* \alpha') =
3(\alpha,\alpha') \mbox{ for all } \alpha,\alpha' \in L_A.$$ 

If $\alpha \in L_A^{\perp} \subset H_2(A,\Z) \subset H_2(\tilde{A},\Z)$,
then $q_*(\alpha) \in L_X$ with
$$(\xi,q_*(\alpha))=(q^*(\xi),\alpha)=0$$
for all $\xi \in L_X$. This implies $\tilde{q}_*(\alpha)=0$.
 
And conversely, if $\tilde{q}_*(\alpha)=0$, then $(\tilde{q}^*\xi,\alpha)=
(\xi,\tilde{q}_*(\alpha))=0$ for all $\xi \in L_X$. So
$\alpha \in L_A^{\perp}$. 
Hence
$\tilde{q}_*$ defines an injective map
$$ \tilde{q}_*:H_2(A,\Z)/L_A^{\perp} \simeq \check{L}_A \to L_X.$$

The discriminant of the image lattice is
$$3^4 \cdot d(\check{L}_A) = -3^3 = d(L_X).$$
This shows $\tilde{q}_*(\check{L}_A)=L_X.$ \qed

{\bf Corollary:} {\em The lattice $L_X$ has an integral
basis, in which its quadratic form has the matrix}
$$\left( \begin{array}{rrrr}
 0 & 3 & 0 & 0 \\ 3 & 0 & 0 & 0 \\ 
 0 & 0 & 6 & 3 \\ 0 & 0 & 3 & 2 \\ \end{array} \right).$$

\VV
\section{Moduli}

\subsection{The period on $A$}

A complex structure on the real four-dimensional torus 
$A=H_1(A,\R)/H_1(A,\Z)$ is given by two complex coordinates $z_1,z_2$.
These are $\R$-linear maps $z_k:H_1(A,\R) \to \C$ inducing an
$\R$-isomorphism $(z_1,z_2):H_1(A,\R) \to \C^2$. This complex
structure determines a {\em period} 
$\omega=z_1 \wedge z_2 \in H^2(A,\C)$. $\omega$ is uniquely
determined by the complex structure up to multiplication
by complex scalars. This period satisfies the
following two {\em period relations}
$$(i) \quad \omega \wedge \omega =0, \qquad 
(ii) \quad \omega \wedge \bar{\omega} > 0.$$

Relation $(i)$ is obvious. Relation $(ii)$ means the following:
Write $z_k=x_k+i y_k$ with $\R$-linear maps 
$x_k,y_k:H_1(A,\R) \to \R$. Then
\begin{eqnarray*}
\omega \wedge \bar{\omega}
 &=& (x_1+iy_1)\wedge (x_2+iy_2) \wedge (x_1-iy_1) \wedge (x_2-i y_2) \\
 &=& -(x_1+i y_1)\wedge (x_1-i y_1)\wedge (x_2+i y_2)\wedge (x_2-i y_2) \\
 &=& -[ -2i \cdot (x_1 \wedge y_1) \wedge (-2i) \cdot (x_2 \wedge y_2) ] \\
 &=& 4 \cdot x_1 \wedge y_1 \wedge x_2 \wedge y_2
\end{eqnarray*}
is a positive multiple of the form  
$ x_1 \wedge y_1 \wedge x_2 \wedge y_2$ defining the orientation
by the complex structure. 

There is also the converse:

{\em Given a class $\omega \in H^2(A,\C)$ satisfying the
period relations $(i)$ and $(ii)$, there is a unique complex
structure on $A$ belonging to this form.}

Proof. Since $\omega \wedge \omega = 0$ by $(i)$, the form
$\omega \in \Lambda^2 H^1(A,\C)$ decomposes, say
$\omega = z_1 \wedge z_2$ with 
$z_1,z_2 \in H^1(A,\C)=Hom_{\R}(H_1(A,\R),\C)$. 
These functions $z_1$ and $z_2$ are uniquely determined by $\omega$
up to complex linear combination.
Write
$z_k=x_k+i y_k$ as above. Then $(ii)$ implies 
$x_1 \wedge y_1 \wedge x_2 \wedge y_2 \not= 0$ and the
map $(z_1,z_2):H_1(A,\R) \to C^2$ is bijective, i.e., it
defines a complex structure on $A$. \qed

\vv
{\bf Proposition.} {\em If $t^* \omega = \omega$, then the
map $t:A \to A$ is $\C$-linear.}

Proof. By assumption
$$t^*(z_1) \wedge t^*(z_2) = t^*(\omega) = \omega = z_1 \wedge z_2.$$
So $t^*(z_1)$ and $t^*(z_2)$ are complex linear
combinations of $z_1$ and $z_2$. \qed

{\bf Corollary.} {\em Isomorphism classes $(A,t)$ of complex tori $A$
with an order-three automorphism $t$ as in 1.1 are classified by 
the period domain $\Omega/G$ where

$\Omega = \{ \C \cdot \omega \in \P(L_A): \;
\omega \wedge \omega=0, \, \omega \wedge \bar{\omega} >0\}$, 

$G$ is the group of orientation-preserving $\Z$-isomorphisms 
$g:L_A \to L_A$ commuting with the $t$-action
specified in 1.1.}

Proof. Let $(A_1,t_1)$ and $(A_2,t_2)$ be two such complex tori with
automorphisms. Fix isomorphisms $H_1(A_k,\Z)=\Z^4$  such that $t_k$
acts as in 1.1. An isomorphism $\phi:A_1 \to A_2$ preserves the
automorphism if
$$\phi \circ t_1 = t_2 \circ \phi.$$
It induces an isomorphism $\phi_*:\Z^4 \to \Z^4$ commuting with $t$. 
So $\phi_* \in G$. It sends $\Omega$ to $\Omega$ and
$\omega_1$ to $\omega_2$.   \qed

\vv
This group $G$ can be descibed a little more explicitly: Indeed
$$g = \left( \begin{array}{cc} A & B \\ C & D \\ \end{array} \right)$$
with $A,B,C,D$ integral $2 \times 2$-matrices belongs to
$G$, if it is invertible, of determinant $1$, and satisfies
$$AT=TA,...,DT=TD, \mbox{ with } T = \left( \begin{array}{rr}
                   0 & -1 \\ 1 & -1 \\ \end{array} \right).$$
And the set of $2 \times 2$-matrices commuting with $T$ consists
of the $\Z$-algebra generated by the unit matrix $\E$ and $T$.

\vv
{\bf Theorem.} {\em The moduli space $\Omega/G$ is connected,
of dimension two.}

Proof. It suffices to show that $\Omega$ is connected. The
intersection form on $L_A$ has signature $(3,1)$. Let us choose real
coordinates $x_1,...,x_4$ on $L_A \otimes \R$ such that
in these coordinates this form is
$$ x_1^2+x_2^2-x_3x_4.$$
Let $z_k=x_k+i \cdot y_k$ be the corresponding
complex coordinates on $L_A \otimes \C = H^2(A,\C)^{inv}$.
If $\omega$ has the coordinates $(c_1,...,c_4)$, then
$$\omega \wedge \omega = c_1^2+c_2^2-c_3c_4=0.$$
Now, if $c_3=0$, then $c_1^2+c_2^2=0$ too. The second period
condition in this case is
$$\omega \wedge \bar{\omega}=|c_1|^2+|c_2|^2 >0,$$
satisfied unless $\omega=(0,0,0,1)$.

If $c_3 \not=0$ we may assume $c_3=1$. Hence 
$\omega=(z_1,z_2,1,z_1^2+z_2^2)$ and
\begin{eqnarray*}
\omega \wedge \bar{\omega} &=& |c_1|^2+|c_2|^2 -\Re (c_1^2+c_2^2) \\
   &=& 2(\Im(c_1)^2+\Im(c_2)^2) \\
   &>& 0
\end{eqnarray*}
unless $\Im(c_1)=\Im(c_2)=0$. So $\Omega \cap \{c_3 \not=0\}$
is just a copy of $\C^2 \setminus \R^2$, hence connected.
And $\Omega \cap \{c_3=0\}$ lies in its boundary. \qed

\vv
The condition $\omega \wedge \omega=0$ defines a non-degenerate
quadric in $\P_3=\P(L_A \otimes \C)$. It is not contained in
any hyper-plane. So the open subset $\Omega$ of this quadric
is not contained in a hyper-plane too. 

Now let $A$ (and $X$) be algebraic with $C \subset A$ the
pull-back of some ample divisor on $X$. It determines
a class $\gamma_C \in L_A$. Since 
$$\gamma_C \wedge \omega \sim \int_C \omega = 0,$$
in this case the period $\omega$ lies in the hyper-plane
$\gamma_C^{\perp} \subset L_A \otimes \C$. This proves

{\bf Theorem.} {\em The general complex torus $A$ with an order-three
automorphism $t$ as in 1.1 is not algebraic.}

\vv
There is also this converse: If a class $\gamma \in H^2(A,\Z)$ has the 
property $\gamma \wedge \omega=0$, then it is of type $(1,1)$. 
Being integral it is the first chern class of a line bundle
on $A$. For $\omega \in L_A \otimes \C$, in particular all
classes $\gamma$ in the two-dimensional lattice
$L_A^{\perp} \subset H^2(A,\Z)$ have this property.
This proves

{\bf Proposition.} {\em For all complex tori $A$ with an
automorphism $t$ as above, the group
$$Pic(A)/Pic^0(A)$$
has rank at least two.}

But beware: The lattice $L_A^{\perp}$ is negative definite.
So all classes $\gamma \in L_A^{\perp}$ belong to line
bundles on $A$, but if $A$ is not algebraic, these
line bundles do not come from divisors.

\vV 
\subsection{Polarizations on $X$}

Let $\xi_1,...,\xi_4 \in L_X \subset H^2(\tilde{X},\Z)$ the
basis from the corollary in 1.3. Each class $n_1\xi_1+...+n_4 \xi_4,
\, n_i \in \Z,$
has self-intersection
$$6 n_1n_2 + 6 (n_3^2+n_3n_4) + 2 n_4^2   
\equiv 0 \mbox{ or } 2 \; mod \; 6.$$

{\bf Theorem.} {\em There is no quartic surface
$X \subset \P_3$ with nine cusps and no other
singularities.}

Proof. A general hyperplane section $C \subset X$
would define a divisor class $\xi=[C] \in J_X$ with 
self-intersection $4 \not\equiv 0 \mbox{ or }2 \; mod \; 6$. \qed 

\vv
On the other hand, for each integer $d \equiv 0 \mbox{ or } 2 \; mod \; 6$
there are $K3$-surfaces with nine cusps carrying a 
polarization of degree $d$: If $d \equiv 0 \; mod \; 6$, then put
$n_1=1, \, n_3=n_4=0$ and $n_2=d/6$.  
If $d=6k+2$, then
put $n_1=n_4=1, \, n_3=0$ and $n_2=(d-2)/6$. 
In fact, [BL] explicitly describe all surfaces with polarization of
degree two, and [BV] gave explicit examples of surfaces with 
polarizations of degree six.  

\vV
\subsection{Polarizations on $A$}

Fix some class $\alpha \in L_A$ with $\alpha \wedge \alpha>0$.
It defines a hyper-plane $\alpha^{\perp} \subset L_A \otimes \C$
and a curve $\Omega_{\alpha}=\Omega \cap \alpha^{\perp}$.
For all periods $\omega \in \Omega_{\alpha}$ the complex
torus $A$ defined by $\omega$ carries line-bundles 
of class $\alpha$. Since $\alpha \wedge \alpha>0$, these
line-bundles are ample, and $A$ is algebraic.

{\bf Proposition.} {\em All these complex curves 
$\Omega_{\alpha} \subset \Omega$ consist of two
connected components. In fact, they are a union of two copies
of the complex upper half-plane.}

Proof. As $\alpha \wedge \alpha > 0$, the lattice
$L_{\alpha}:=L_A \cap \alpha^{\perp}$ has signature $(2,1)$. So there are 
real coordinates $x_1,x_2,x_3$ of $L_{\alpha} \otimes \R$
in which the intersection form is 
$x_1^2-x_2x_3$.
Let $z_k$ be the corresponding complex coordinates on 
$L_{\alpha} \otimes \C$. 
Each $\omega \in \Omega_{\alpha}$ has coordinates 
$$(c_1,c_2,c_3), \quad  c_k \in \C$$
with
$$ \omega \wedge \omega = c_1^2 -c_2c_3 = 0.$$
If $c_3=0$, then $c_1=0$ implies $\omega \wedge \bar{\omega}=0$,
a contradiction. So $c_3 \not=0$. We may assume $c_3=1$
and $\omega=(c_1,c_1^2,1)$.

Then
$$\omega \wedge \bar{\omega} = |c_1|^2 -\Re(c_1^2) =
2 \Im(c_1)^2 >0$$
unless $\Im(c_1)=0$. It follows that
$\Omega_{\alpha}=\C \setminus \R$. \qed
 
\vV
\subsection{The classification problem}

Of course, it would be nice to have a classification for
the set of abelian surfaces with an automorphism $t$ of order
three as in section 1.1, together with a $t$-invariant polarization.
Polarizations of degree $2d$ on abelian surfaces 
(without $t$) are classified by their
elementary divisors 
$$d_1,d_2 \in \N, \quad d_1|d_2, \quad 2 \cdot d_1d_2 = 2d$$
in the sense that they belong to the same connected moduli
space, if these elementary divisors coincide. However for fixed
$d_1,d_2$, not all pairs $(A,t)$ admitting a 
polarization of type $(d_1,d_2)$ are toplogically equivalent. 
There is an obvious reason: If two polarizations, i.e. primitive
vectors $\alpha$ and $\alpha' \in L_A$ with $\alpha^2=(\alpha')^2 >0$
are toplogically equivalent, they are conjugate under the
arithmetic group $G$ from 2.1. In particular they must be conjugate 
under the orthogonal group $O(L_A)$. And then 
$\alpha \in L_A \subset \check{L}_A$ is primitive, if and only if 
$\alpha'$ is prinitive in $\check{L}_A$.

It is easy to see that the vector 
$$\alpha:=\gamma_1+3\gamma_2$$ 
of length $\alpha^2=6$ is primitive in $\check{L}_A$, while the vector
$$\alpha'=\gamma_3+\gamma_4$$
of the same length $6$ is not. So the polarizations $\alpha$ and $\alpha'$
are not topologically equivalent. (In fact, by 1.3, the polarization
$\alpha'$ descends to the quotient $X=A/t$, while the 
$t$-invariant  polarization $\alpha$ does not do this).

The real classification problem is this:

{\bf Problem:} {\em Classify the $G$-orbits on the sets of primitive
vectors $\alpha \in L_A$ of the same length $\alpha^2>0$ and
primitive / not primitive in $\check{L}_A$.}

I expect the groups $G$ and $O(L_A)$ to be more or less the same, but
I have no idea, whether they act transitively or not on primitve
vectors $\alpha \in L_A$ of the same length, primitive, resp.
not primitive in $\check{L}_A$. 

\vv
Each abelian surface ($=$ algebraic torus) 
$A$ with automorphism $t$ has Picard number
three. 
Let me give examples
of simple abelian surfaces $A$ with automorphism $t$

To do this, I have to identify the orthogonal complement
$L_A^{\perp} \subset H^2(A,\Z)$:

{\bf Proposition.} {\em The lattice $L_A^{\perp}$ admits
a $\Z$-basis $\delta_1,\delta_2$ in which the quadratic form is given
by the matrix}
$$\left(\begin{array}{rr} -2 & 1 \\ 1 & -2 \\ \end{array} \right).$$ 

Proof. Recall the basis $\alpha_1,\beta_1,\alpha_2,\beta_2$
from section 1.1. Put 
$$\delta_1:=\alpha_1\wedge \alpha_2- \beta_1 \wedge \beta_2,
\quad \delta_2=\alpha_1\wedge \beta_2+\beta_1 \wedge \alpha_2+
\beta_1 \wedge \beta_2.$$
One easily checks that these classes are orthogonal to the
basis $\gamma_1,...,\gamma_4$ of $L_A$ from section 1.2, and
$$\delta_1\wedge \delta_1 = \delta_2 \wedge \delta_2 = -2,
\; \delta_1 \wedge \delta_2 = 1.$$
So $\delta_1$ and $\delta_2$ span a sublattice of $L_A^{\perp}$
in which the form has the matrix above.
Its discriminant is $3$. This sublattice therefore is primitive, and must
coincide with $L_A^{\perp}$. \qed

\vv
The essential observation is: This lattice does not represent 
$-12 \cdot n^2$ for any $n \in \N$.

Proof. Assume there are $k,l \in \Z$ with 
$$(k \delta_1+l \delta_2)^2=-2(k^2+l^2-kl) =-12n^2.$$
We solve the quadratic equation
$$k^2-kl+l^2-6n^2=0$$
for $k$ to find the two solutions
$$k_{1,2}=\frac{1}{2}(l\pm \sqrt{24n^2-3l^2}).$$
This shows 
$$24n^2-3l^2=(2k_{1,2}+l)^2=w^2$$
is a square. Obviously $3$ divides $w$, so write
$w=3w'$ and
$$8n^2-l^2=3(w')^2.$$
If neither $n$ nor $l$ are divisible by three,
we find
$$8n^2 \equiv 2, \quad l^2 \equiv 1  \quad modulo \; 3,$$
impossible. So $n=3n'$ and $l=3 l'$. But this leads to
$$24(n')^2-3(l')^2 = (w')^2,$$
the original equation. As we may repeat this argument infinitely often,
this is a contradiction. \qed

{\bf Proposition.} {\em There are abelian surfaces with an automorphism $t$ 
as above, and not carrying elliptic curves.}

Proof. Recall the basis $\gamma_1,...,\gamma_4 \in L_A$
from section 1.2. Clearly 
$$\gamma:=-\gamma_1+6 \gamma_2+2(\gamma_3+\gamma_4)$$
is a primitive vector of length $12$. Let $A$ be a surface
carrying a polarization with class $\gamma$, and with
Neron-Severi group of rank three. Then this 
Neron-Severi group is spanned by $\gamma, \delta_1$ and
$\delta_2$. If it would carry elliptic curves,
there would be a class $n \gamma + k \delta_1+l \delta_2$
of length
$$12 n^2 - 2(k^2-kl+l^2)=0,$$
in conflict with the observation above. \qed

\VV
\section{Appendix: Quartic surfaces with eight cusps}

In this section I compute the equations of those
quartic surfaces in $\P_3(\C)$, which are the projections
of the sextic surfaces in $\P_4$ from [BV] out of one of their
nine cusps. Then I show by direct computation,
that these projected surfaces have precisely eight cusps, and no further
singularities.

So recall the sextic surfaces from [BV]. On p.141 they are presented 
as a complete intersection of a hyper plane $P$, a quadric $Q$ and
a cubic $C \subset \P_5$ with
$$ \begin{array}{ll}
P: & x_1+x_2+x_3+x_4+x_6+x_6=0,\\
Q: & (1+k)(x_1x_2+x_1x_3+x_2x_3)+(1-k)(x_4x_5+x_4x_6+x_5x_6)=0, \\
C: & (1+k)^2x_1x_2x_3+(1-k)^2x_4x_5x_6=0.\\
\end{array}$$
Here $k \not=0$ is a complex parameter. The nine cusps of these
surfaces are the points
$$(1:0:0-1:0:0) \quad \mbox{etc.}$$
with one coordinate $x_1,x_2$ or $x_3$ equal to $1$ and 
one coordinate $x_4,x_5$ or $x_6$ equal to $-1$. 

I fix the cusp $(1:0:0:-1:0:0)$ and project from it onto the
$3$-plane 
$$x_1=x_4=-\frac{1}{2}(x_2+x_3+x_5+x_6)=:\sigma.$$
the rays of projection are parametrized by
$$(\lambda+\mu s: \mu x_2:\mu x_3:-\lambda+\mu \sigma: \mu x_5: \mu x_6)
\quad \mbox{with} \quad (\lambda:\mu) \in \P_1.$$
On this ray the equation of $Q$ restricts to
$$\lambda[(1+k)(x_2+x_3)-(1-k)(x_5+x_6)]
+\mu[(1+k)(\sigma(x_2+x_3)+x_2x_3)+(1-k)(\sigma(x_5+x_6)+x_5x_6)]=0$$
and the equation of $C$ becomes
$$\lambda[(1+k)^2x_2x_3-(1-k)^2x_5x_6]+
\mu \cdot \sigma [(1+k)^2x_2x_3+(1-k)^2x_5x_6]=0.$$
Eliminating $(\lambda:\mu)$ from these two equations,
and replacing the coordinates $x_2,x_3,x_4,x_5$
by $x_0,x_1,x_2,x_3$ leads to the equation
\begin{eqnarray*}
(1+k)^3x_0^2x_1^2+2k(1-k^2)x_0x_1x_2x_3-(1-k)^3x_2^2x_3^2 & & \\
+(1-k^2)(x_0+x_1+x_2+x_3)[(1-k)x_2x_3(x_0+x_1)-
                          (1+k)x_0x_1(x_2+x_3)] &=&0
\end{eqnarray*}
for the projected surface.

The eight nodes of the sextic surface, different from
the center of projection, go onto the four coordinate
vertices
$$(1:0:0:0),...,(0:0:0:1)$$
and the four points 
$$(1:0:-1:0),\,(1:0:0:-1),\,(0:1:-1:0),\,(0:1:0:-1)$$
on the coordinate lines in the plane $x_0+x_1+x_2+x_3=0$.

The center of projection blows up to the pair of lines
$$ \begin{array}{lclcl}
(1+k)x_0-(1-k)x_2 &=& (1+k)x_1-(1-k)x_3 &=& 0 \\
(1+k)x_0-(1-k)x_3 &=& (1+k)x_1-(1-k)x_2 &=& 0 \\
\end{array}$$
on the projected surface.

The projected surfaces have the obvious 
$\Z_2 \times \Z_2$-symmetry
$$x_0 \leftrightarrow x_1 \quad \mbox{and} \quad 
x_2 \leftrightarrow x_3,$$
while the symmetry
$$(x_0,x_1) \leftrightarrow (x_2,x_3), \quad 
k \leftrightarrow (-k)$$
interchanges two surfaces in the family (if $k \not=0$).

\vv
{\bf Claim.} {\em For $k \not=\pm 1$ the quartic surfaces 
have no other singularities than the eight ones specified.}

To have some (not very big) computational advantages,
let me pass to the coordinates
$$y_0=(1+k)x_0,\, y_1=(1+k)x_1, \; y_2=(1-k)x_2, \, y_3=(1-k)x_3,$$
in which the eqution (multiplied by $(1-k^2)$) takes the form
\begin{eqnarray*}
f(y) &=& (1-k)y_0^2y_1^2+2ky_0y_1y_2y_3-(1+k)y_2^2y_3^2+ \\
     & & ((1-k)(y_0+y_1)+(1+k)(y_2+y_3))
         ((y_0+y_1)y_2y_3-(y_2+y_3)y_0y_1) \\
     &=& 0.
\end{eqnarray*}
We have to compute the derivative  
\begin{eqnarray*}
\partial_0 f &=& 2(1-k)y_0y_1^2+2ky_1y_2y_3
                 +(1-k)((y_0+y_1)y_2y_3-(y_2+y_3)y_0y_1)+ \\
             & & ((1-k)(y_0+y_1)+(1+k)(y_2+y_3))(y_2y_3-(y_2+y_3)y_1)
\end{eqnarray*}
To take advantage of the symmetries, let me abbreviate
$$s:=y_0+y_1, \, t:=y_2+y_3, \; p:=y_0y_1, \, q:=y_2y_3.$$
Then
$$\partial_0f = 2((1-k)p+kq)y_1+(1-k)(sq-tp)+
                ((1-k)s+(1+k)t)(q-ty_1),$$
and by the symmetries 
$$\partial_1f = 2((1-k)p+kq)y_0+(1-k)(sq-tp)+
                ((1-k)s+(1+k)t)(q-ty_0).$$
The difference of these two derivatives is
$$\partial_0f-\partial_1f = (y_1-y_0) \cdot
(2(1-k)p+2kq-((1-k)s+(1+k)t)t).$$

So, in a singularity, either
$$y_0=y_1,$$
or
$$((1-k)s+(1+k)t)t = 2(1-k)p+2kq.$$
Inserting the latter in the expression for $\partial_0f$
leads to
$$\partial_0f = 2(1-k)sq-(1-k)tp+(1+k)tq.$$
Using the symmetry $(x_0,x_1,k) \to (x_2,x_3,-k)$ we find: 
At a singular point 
$$
\begin{array}{llrcl}
\mbox{either }x_0=x_1 & \mbox{or} &
    -(1-k)t \cdot p +(2(1-k)s+(1+k)t) \cdot q &=& 0, \\ 
                      &           &
   2(1-k)\cdot p +2k \cdot q &=& ((1-k)s+(1+k)t)t, \\
\mbox{either }x_2=x_3 & \mbox{or} &
    (2(1+k)t+(1-k)s) \cdot p -(1+k)s \cdot q &=& 0, \\
                      &           &
    -2k \cdot p +2(1+k) \cdot q &=& ((1+k)t+(1-k)s)s. \\
\end{array}$$
The two homogeneous equations for $p$ and $q$ form a system with determinant
$$D= (1-k^2)st -(2(1-k)s+(1+k)t)(2(1+k)t+(1-k)s)
=-2((1-k)s+(1+k)t)^2.$$
So, if they hold, either
$$\frac{s}{1+k}+\frac{t}{1-k}=x_0+x_1+x_2+x_3=0,$$
or $p=q=0$. In this case the two inhomogeneous equations for $p$ and $q$ show
$$((1-k)s+(1+k)t)\cdot t = ((1-k)s+(1+k)t) \cdot s = 0.$$
And this again implies $x_0+x_1+x_2+x_3=0$.

So, if a surface for $k \not= \pm1$ is singular at $(x_0:x_1:x_2:x_3)$, then
we are in one of the following cases:

{\em Case I:} $x_0+x_1+x_2+x_3=0$. But then
$$\begin{array}{lclcl}
\D \frac{\partial_0f-\partial_1f}{2(y_1-y_0)} &=& (1-k)p+kq &=& 0, \\
 & & \\
\D \frac{\partial_2f-\partial_3f}{2(y_3-y_2)} &=& -kp+(1+k)q &=& 0. \\
\end{array}$$
This system for $p$ and $q$ has determinant 
$$(1-k^2)+k^2=1,$$
showing $p=q=0$. So $x$ is one of the four points 
$$(1:0:-1:0),\,(1:0:0:-1),\,(0:1:-1:0),\,(0:1:0:-1).$$

{\em Case II:} $x_0=x_1$, hence $y_0=y_1=:y$, but $x_2 \not=x_3$. 
Then there still are the two equations
\begin{eqnarray*}
    (2(1+k)t+2(1-k)y) \cdot y^2 -2(1+k)y \cdot q &=& 0, \\
    -2k \cdot y^2 +2(1+k) \cdot q &=& ((1+k)t+(1-k)2y) \cdot 2y 
\end{eqnarray*}
from the four equations above. Unless $y=0$ this leads to
the two quadratic equations
\begin{eqnarray*}
(1-k)\cdot y^2 +(1+k)t \cdot y -(1+k)q &=& 0, \\
(k-2)\cdot y^2 -(1+k)t \cdot y +(1+k)q &=& 0.
\end{eqnarray*}
Adding both equations we see $y=0$. And if $y=0$, the second equation
shows $q=0$. The point in question is one of the four coordinate
vertices.

{\em Case III:} $x_2=x_3$, but $x_0 \not= x_1$. 
This leads to a coordinate vertex, just like
case II.

{\em Case IV:} $x_0=x_1$ and $x_2=x_3$, hence
$y_0=y_1=:y$ and $y_2=y_3=:z$. Then
$$\partial_0f = 2(1-k)y^3+2k yz^2+(1-k)(2yz^2-2y^2z)
                +2((1-k)y+(1+k)z)(z^2-2yz) = 0,$$
and
$$ \begin{array}{rcccl}
\D \frac{1}{2} \partial_0f &=& (1-k)y^3-3(1-k)y^2z-3kyz^2+(1+k)z^3 &=& 0, \\
 & & & & \\
\D \frac{1}{2} \partial_2f &=& (1-k)y^3+3ky^2z-3(1+k)yz^2+(1+k)z^3 &=& 0. \\
\end{array}$$
The difference of both these equations is
$$3yz \cdot (z-y) =0.$$
Now $yz=0$ leads to the contradiction $y=z=0$. And for $y=z$ both equations
become $-y^3=0$, again a contradiction.

This proves the claim. \qed

\vv 
Next I want to show that the eight singularities indeed are cusps.
To recognize them as cusps, I use the recognition principle from
[BW]. In proper coordinates $u_1,u_2,u_3$ one has to assign weights
$$wt(u_1)=\frac{1}{3}, \; wt(u_2)=wt(u_3)=\frac{1}{2}$$ 
to these coordinates, and to show:

\begin{itemize}
\item The Taylor expansion of the equation of the surface at the 
a singularity contains only monomials of weight $\geq 1$;
\item the weight-$1$ part of the expansion is $u_1^3+u_2u_3$.
\end{itemize}

Let me call this recognition principle SQH-criterion.

\vv
{\bf Claim:} {\em For $k \not= \pm 1$ the coordinate vertices are
cusps on our surfaces.}

By the symmetries, it suffices to consider the point $(1:0:0:0)$.
We consider
\begin{eqnarray*}
f(1,y_1,y_2,y_3) &=& (1-k)y_1^2+2ky_1y_2y_3-(1+k)y_2^2y_3^2+ \\
 & & ((1-k)(1+y_1)+(1+k)((y_2+y_3))
     ((1+y_1)y_2y_3-(y_2+y_3)y_1) \\
 &=& (1-k)(y_1^2+y_2y_3-y_1(y_2+y_3))+ \\
 & & \mbox{terms of order } \geq 3.
\end{eqnarray*}
So the second-order term splits as
$$(k-1)(y_1-y_2)(y_1-y_3).$$
Now we substitute
$$u_1:=y_2+y_3, \; u_2:=y_1-y_2, \; u_3:=y_1-y_3.$$
assigning to the
new variables the weights
$$wt(u_1)=\frac{1}{3}, \; wt(u_2)=wt(u_3)=\frac{1}{2}.$$
We have to collect all terms in the Taylor expansion of $f$ of weight
$\leq 1$. Since all coordinates $y_k$ have weight $\geq 1/3$, 
all fourth order monomials in $y_1,y_2,y_3$ will have weight $>1$.
It suffices therefore to consider the third order term
\begin{eqnarray*}
 & & 2ky_1y_2y_3+(1-k)y_1y_2y_3+(1-k)y_1(y_2y_3-(y_2+y_3)y_1)+ \\
 & & (1+k)(y_2+y_3)(y_2y_3-(y_2+y_3)y_1) \\
 &=& 2y_1y_2y_3-(1-k)y_1^2(y_2+y_3)-(1+k)y_1(y_2+y_3)^2
     +(1+k)(y_2+y_3)y_2y_3 \\
 &=& (\D \frac{1}{4}-\frac{1-k}{4}-\frac{1+k}{2}+\frac{1+k}{4})u_1^3 \\
 & & + \mbox{terms of weight }>1 \\
 &=& \D -\frac{1}{4}u_1^3 \\
 & \not=& 0.
\end{eqnarray*}
The term of weight $1$ in the Taylor-expansion is
$$(k-1)u_2u_3-\frac{1}{4}u_1^3.$$
By [BW], indeed the singularity is a cusp.

\vv
{\bf Claim:} {\em For $k \not= \pm 1$ the other four singularities
have type $A_2$ too.}

By the symmetries, it suffices to consider the singularity
$(1:0:-1:0)$ only. We use the coordinates $x_0,y_1,x_2,y_3$ to
write
\begin{eqnarray*}
f(x_0,y_1,x_2,y_3)
 &=& (1+k)x_0^2y_1^2+2kx_0x_2y_1y_3-(1-k)x_2^2y_3^2+ \\
 & & ((1-k^2)(x_0+x_2)+(1-k)y_1+(1+k)y_3) \cdot \\
 & & (x_0x_2(y_3-y_1)+x_2\frac{y_1y_3}{1+k}-x_0\frac{y_1y_3}{1-k}). 
\end{eqnarray*}
Now we substitute
$$x_0=u+v, \; x_2 = u-v$$
and write $f$ in the form
\begin{eqnarray*}
f(u,v,y_1,y_3)
 &=& (1+k)(u+v)^2y_1^2+2k(u^2-v^2)y_1y_3-(1-k)(u-v)^2y_3^2+ \\
 & & (2(1-k^2)u+(1-k)y_1+(1+k)y_3) \cdot \\
 & & ((u^2-v^2)(y_3-y_1)+(u-v)\frac{y_1y_3}{1+k}-(u+v)\frac{y_1y_3}{1-k}). 
\end{eqnarray*}
The Taylor-expansion at the singularity is
\begin{eqnarray*}
f(u,1,y_1,y_3)
 &=& (1+k)(u+1)^2y_1^2+2k(u^2-1)y_1y_3-(1-k)(u-1)^2y_3^2+ \\
 & & (2(1-k^2)u+(1-k)y_1+(1+k)y_3) \cdot \\
 & & ((u^2-1)(y_3-y_1)+(u-1)\frac{y_1y_3}{1+k}-(u+1)\frac{y_1y_3}{1-k}).\\
 &=& (1+k)y_1^2-2ky_1y_3-(1-k)y_3^2+ \\
 & & (2(1-k^2)u+y_1+y_3-k(y_1-y_3)) \cdot (y_1-y_3)+ \\ 
 & & \mbox{terms of order } \geq2 \\
 &=& y_1^2-y_3^2 +(2(1-k^2)u+y_1+y_3)(y_1-y_3) \\
 &=& 2(y_1-y_3)((1-k^2)u+y_1+y_2).
\end{eqnarray*}
We substitute
$$z_1:=(1-k^2)u+y_1+y_3, \; z_2=y_1-y_3$$
and consider the third-order term of the Taylor expansion
$$2(1+k)uy_1^2+2(1-k)uy_3^2 +(z_1+(1-k^2)u-kz_2)\cdot 
  (-\frac{2}{1-k^2}y_1y_3).$$
We assign to the variables $z_1,z_2,u$ the following weights
$$wt(z_1)=wt(z_2)=\frac{1}{2}, \; wt(u)=\frac{1}{3}.$$
Since each monomial of degree four then has weight $>1$, it suffices
indeed to treat the third order term only.
Its part of weight $1$ is
\begin{eqnarray*}
 & &2u((1+k)\frac{(1-k^2)^2}{4}u^2+(1-k)\frac{(1-k^2)^2}{4}u^2)
-(1-k^2)u\cdot \frac{2}{1-k^2} \cdot \frac{1}{4}(1-k^2)^2u^2 \\
 &=& 2u \cdot \frac{(1-k^2)^2}{2} u^2 - u \frac{(1-k^2)^2}{2}u^2 \\
 &=& \frac{1}{2}(1-k^2)^2 u^3 \\
 &\not=& 0.
\end{eqnarray*}
The weight-$1$ part of $f$ is
$$z_1z_2-\frac{1}{4}(1-k^2)^2u^2,$$
and again the SQH-criterion [BW] shows that the singularity is a cusp.

\vV
\section{References}

\noindent
[B] Barth, W.: K3-surfaces with nine cusps, to appear
in Geom. Dedic.

\noindent
[BL] Birkenhake C., Lange, H.: A family of abelian surfaces and
curves of genus four. manuscr. math. 85, 393-407 (1994)  

\noindent
[BV] Bertin, J., Vanhaecke, P.: The even master system
and generalized Kummer surfaces. Math. Proc. Camb. Phil. Soc.
116, 131-142 (1994)

\noindent
[BW] Bruce, J.W., Wall C.T.C.: On the classification of
cubic surfaces. J. London Math. Soc. 19, 245-249 (1979)

\noindent
[N1] Nikulin, V.V.: On Kummer surfaces. Math. USSR Izv. 9, No 2, 261-275
(1975)

\noindent
[N2] Nikulin, V.V: Integral symmetric bilinear forms and some of
their applications. Mat. Ussr. Izv. 14, No 1, 103-167 (1980)

\noindent
[U] Urabe, T.: Elementary transformations of Dynkin graphs
and singularities on quartic surfaces. Invent. math. 87, 549-572
(1987) 
 
\end{document}